\documentclass[letterpaper, 10 pt, conference]{ieeeconf}  

\IEEEoverridecommandlockouts                              
\usepackage{cite}
\usepackage{float}
\usepackage{multicol}
\usepackage{multirow}
\usepackage{hyperref}
\usepackage{graphicx}
\usepackage{amsfonts}
\usepackage{amsmath}
\usepackage{color}
\usepackage{epstopdf}
\usepackage{subcaption}
\usepackage{multirow}
\usepackage{bbm}
\usepackage{bm}
\usepackage{comment}

\title{\LARGE \bf
Chance Constrained Optimal Power Flow with Primary Frequency Response 
}

\author{ \parbox{3 in}{\centering Michael Chertkov \\
        Theory Division \& Center for Nonlinear Studies\\
        Los Alamos National Laboratory\\
        Los Alamos, NM 87545  \\
         {\tt\small chertkov@lanl.gov}}
         \hspace*{ 0.5 in}
         \parbox{4 in}{ \centering Yury Dvorkin \\
        Dep. of Electrical \& Computer Engineering \\
         Tandon School of Engineering\\
         New York University\\
         New York, NY 11201 \\
         {\tt\small dvorkin@nyu.edu}}
}

\begin{document}

\maketitle

\thispagestyle{empty}
\pagestyle{empty}

\vspace{-10mm}
\begin{abstract}

Primary frequency response  is  provided by synchronized generators through their speed-droop governor characteristic in response to instant frequency deviations that exceed a certain threshold, also known as the governor dead zone. 
This paper presents an Optimal Power Flow (OPF) formulation that explicitly  models i) the primary frequency response constraints using a nonlinear speed-droop governor characteristic with a dead zone and ii) chance constraints on power outputs of conventional generators and line flows imposed by the uncertainty of renewable generation resources. The proposed formulation is evaluated and compared to the standard CCOPF formulation on a modified 118-bus IEEE Reliability Test System.

\end{abstract}

\section{Introduction} \label{sec:intro}


As a result of policy initiatives and incentives, renewable generation  resources have already exceeded 10\% penetration levels, in terms of annual electricity produced, in some interconnections; even higher targets are expected to be reached in the forthcoming years \cite{malik_2014}. Integrating renewable generation resources increases reserve requirements needed to deal with their uncertainty and variability and, at the same time, tends to replace conventional generation resources that are most suitable for the efficient reserve provision, \cite{troy_2010}. Among these reserve requirements, the ability to provide sufficient \textit{primary frequency response}, i.e.   the automatic governor-operated response of synchronized generators  intended to compensate a sudden power mismatch based on local frequency measurements, is the least studied \cite{eto_2012}.  The impacts of renewable generation resources on the secondary  and tertiary reserve requirements are reviewed in \cite{makarov_2009, dvorkin_2014c}.

Data-driven studies in \cite{du_2014} and \cite{ingleson_2010} manifestly reveal that the primary frequency response in major US interconnections has drastically reduced over the past decades and attribute this effect to increased penetrations of renewable generation resources. Primary frequency response constraints are modeled  in \cite{omalley_2005d} and \cite{restrepo_2005b}. These studies consider a nonlinear speed-droop governor characteristic with an intentional dead zone that makes it possible to preserve the primary frequency response for reacting to relatively large frequency deviations caused by sudden generation and demand failures \cite{jaleeli_1992}. However, the formulations in \cite{omalley_2005d} and \cite{restrepo_2005b} ignore transmission constraints and the uncertainty of renewable generation resources that may lead to overload and capacity scarcity  events when the primary frequency response is deployed.

To deal with the stochastic nature of  renewable generation resources  and their impacts on the secondary and tertiary reserve requirements, the standard optimal power flow (OPF) formulations have been enhanced using chance constrained programming \cite{bienstock_2014},  scenario-based stochastic programming \cite{yuan_2011b} and robust optimization \cite{jabr_2013}. The formulation in \cite{bienstock_2014} postulates  that the uncertainty and variability of wind power generation resources follow a given Gaussian distribution that makes it possible to convert the Chance Constrained OPF (CCOPF) into a second-order cone program, which is then solved using a cutting-plane-based procedure. Relative to \cite{bienstock_2014}, the formulation in \cite{yuan_2011b} is more computationally demanding since it requires computationally expensive scenario sampling, \cite{dvorkin_2014a}. Unlike \cite{bienstock_2014} and \cite{yuan_2011b}, the formulation in \cite{jabr_2013} disregards the likelihood of individual scenarios within a given uncertainty range and tends to yield overly conservative solutions. Based on the CCOPF formulation in \cite{bienstock_2014}, several extensions have been developed. Thus, \cite{lubin_2016} implements a distributionally robust CCOPF formulation that internalizes parameter uncertainty on the Gaussian distribution characterizing wind power forecast errors as explained in \cite{dvorkin_2016f}.  In \cite{roald_2017}, the CCOPF formulation is extended to accommodate corrective control actions.  The formulation in \cite{roald2016c} describes the CCOPF formulation that uses wind curtailments for self-reserve to reduce the secondary and tertiary reserve provision by conventional generation resources. In \cite{roald_2015}, the chance constraints are modified to selectively treat large and small wind power perturbations using weighting  functions. Notably, the convexity of the original formulation in \cite{bienstock_2014} is preserved in \cite{roald_2017, roald2016c, roald_2015, lubin_2016}.

The formulations in \cite{bienstock_2014, roald_2017, roald2016c, roald_2015,lubin_2016} have been proven to reliably and cost-efficiently deal with the uncertainty and variability imposed by renewable generation resources at a computational acceptable cost, even for realistically large networks \cite{bienstock_2014,lubin_2016}. However, these formulations neglect to account for nonlinear primary frequency response policies, i.e. the dead zone of the speed-droop governor characteristic. From the reliability perspective, this may result in the inability to timely arrest a frequency decay caused by credible contingencies, \cite{du_2014, dvorkin_2016, jaleeli_1992}. Furthermore, ignoring nonlinear primary frequency response policies may lead to suboptimal dispatch decisions and cause unnecessary out-of-market corrections that would eventually increase the overall operating cost, \cite{hedman_2015}. This paper {\it proposes a CCOPF formulation with a nonlinear primary frequency response policy that seeks  the least-cost dispatch and  primary frequency response provision  among available conventional generation resources}. The main contributions of this paper are:
\begin{enumerate}
\item The proposed CCOPF-Primary Frequency Rresponse (PFR) formulation enhances the CCOPF formulation in \cite{bienstock_2014} by explicitly considering primary frequency response constraints of conventional generators.
\item As in \cite{omalley_2005d} and \cite{restrepo_2005b}, the  primary frequency response is modeled using a nonlinear speed-droop governor characteristic with an intentional dead zone, i.e. the primary frequency response does not react to frequency fluctuations below a given threshold.
\item The weighted chance constraints from \cite{roald_2015} are used to  reformulate the proposed CCOPF-PFR formulation into a convex program that can be solved using off-the-shelf solvers. The CCOPF-PFR and CCOPF formulations are compared using a modification of the 118-bus IEEE Reliability Test System \cite{ieee118}.
\end{enumerate}


The rest of this paper is organized as follows. Section~\ref{sec:model}  describes the  proposed CCOPF-PFR on the standard CCOPF from \cite{bienstock_2014}. Section~\ref{sec:solution_approach} describes the solution technique. Section~\ref{sec:case_study} compares the CCOPF-PFR and CCOPF formulations quantitatively. Section~\ref{sec:conclusion} summarizes the key findings.

\section{CCOPF-PFR Formulation} \label{sec:model}

This section derives the CCOPF-PFR formulation. Section \ref{subsec:CCOPF} reviews the CCOPF formulation from \cite{bienstock_2014}. Next, Section \ref{subsec:affine} details  a generic affine frequency control policy, which is generalized in Section \ref{subsec:dead} to include a given dead zone. Section \ref{subsec:weighted} applies the weighted chance constraints from \cite{roald_2015} to obtain the final CCOPF-PFR formulation.

\subsection{Standard CCOPF}
\label{subsec:CCOPF}

As customarily done in the analysis of transmission grids, the deterministic OPF based on the DC power flow (PF) approximation is stated as follows:
\begin{eqnarray}
&\min_{p_g,\phi,\theta}& \sum_{i\in{\cal V}_g} C_i(p_i)
\label{DC-OPF}\\
&\mbox{s. t.}& \nonumber\\
&& \sum_{i\in{\cal V}} p_i=0 \label{balance} \\
& \forall i\in{\cal V}_g:& p_i\in [\underline{p}_i,\overline{p}_i] \label{gen_constr_hard}\\
& \forall i\in{\cal V}:& p_i=\sum_{j:\{i,j\}\in{\cal E}} \phi_{ij} \label{flow_node}\\
& \forall \{i,j\}\in{\cal E}: &\phi_{ij}=\beta_{ij}(\theta_i-\theta_j) \label{flow_line}\\
&&\phi_{ij}\in[-\bar{\phi}_{ij},\bar{\phi}_{ij}],\label{line_constr_hard}
\end{eqnarray}
where the transmission grid-graph is ${\cal G}=({\cal V},{\cal E})$, where  ${\cal V}$ and ${\cal E}$ denote the sets of nodes (buses) and edges (transmission lines), and ${\cal V}$ is split into three subsets: conventional generators, ${\cal V}_g$, loads, ${\cal V}_l$, and wind farms, ${\cal V}_w$, i.e. ${\cal V}={\cal V}_g \cup {\cal V}_l \cup {\cal V}_w$. Eq.~(\ref{DC-OPF}) is optimized over the i)  vector of power outputs of conventional generators, $p_g= (p_i|i\in \cal{V}_g)$, ii) vector of transmission power flows, $\phi=(\phi_{ij}=-\phi_{ji}|\{i,j\}\in{\cal E})$, and  iii) vector of voltage angles, $\theta=(\theta_i|i\in {\cal V})$. The input parameters include the minimum, $\underline{p}_i$, and the maximum, $\overline{p}_i$, limits on  the power output of conventional generators, the vector of  nodal loads, $p_l=(p_i|i\in{\cal V}_l)$, the vector of wind power injections, $p_w=\rho=(\rho_i|i\in{\cal V}_w)$, as well as the vector of  line impedances, $(\beta_{ij}, \forall \{i,j\} \in \cal{E})$ and the vector of the line power flow limits, $(\overline{\phi}_{ij},  \forall \{i,j\} \in \cal{E})$. Eq.~\eqref{DC-OPF} minimizes the total operating cost using a convex, quadratic cost function of each conventional generator,  $C_i(\cdot)$. The system-wide power balance is enforced by Eq.~\eqref{balance}.   Eq.~\eqref{gen_constr_hard} limits the power output of the conventional generators. The nodal power balance is enforced by Eq.~\eqref{flow_node} using the line flows computed in Eq.~\eqref{flow_line} based on the {DC} power flow approximation. Eq.\eqref{line_constr_hard}  limits the line flows.

Eqs.~\eqref{DC-OPF}-\eqref{line_constr_hard}  seek the least-cost solution assuming fixed power outputs of wind farms. In fact, these outputs are likely to vary due to the  inherent  wind speed uncertainty and variability \cite{dvorkin_2016f} that can be accounted for as follows by using the chance constrained framework of \cite{bienstock_2014,lubin_2016}:
\begin{eqnarray}
&\min\limits_{p_g^{(0)},\theta, \rho}& \sum_{i\in{\cal V}_g} \mathbb{E}_{\rho}\left[C_i(p_i)\right]
\label{CC-OPF}\\
&\mbox{s. t.}& \nonumber\\
&\forall i\in{\cal V}_g:& \mbox{Prob}_{\rho} \left[p_i\leq \underline{p}_i\right]\leq \varepsilon_i^\downarrow \label{gen_constr_hard1}\\
&& \mbox{Prob}_{\rho} \left[p_i\geq \overline{p}_i\right]\leq \varepsilon_i^\uparrow \label{gen_constr_hard2}\\
& \forall \{i,j\}\in{\cal E}:& \mbox{Prob}_{\rho} \left[\phi_{ij}\leq -\underline{\phi}_{ij}\right]\leq \varepsilon_{ij}^\downarrow\label{line_constr_hard1}\\
&& \mbox{Prob}_{\rho} \left[\phi_{ij}\geq \overline{\phi}_{ij}\right]\leq \varepsilon_{ij}^\uparrow. \label{line_constr_hard2}
\end{eqnarray}
where we follow \cite{dvorkin_2016f} to describe $\rho=p_w$ as a random quantity represented through the exogenous statistics of the wind power. Consistently with \cite{dvorkin_2016f} we use the term ``uncertainty" in this paper to refer to wind power forecast errors (associated with the imperfection of measurement/computation/predicition tools) and the term ``variability" to refer to random fluctuations of wind power output caused by natural atmospheric processes. This statistics is assumed to be Gaussian 
\begin{eqnarray}
& \forall i\in{\cal V}_w: & \mathbb{E}\left[\rho_i\right]=\bar{\rho}_i,\label{wind_mean}\\
&& \ \mathbb{E}\left[\left(\rho_i-\bar{\rho}_i\right)^2\right]=R_i,
\label{wind_cov}
\end{eqnarray}
where $\bar{\rho}_i$ and $R_i$ are the mean and covariance. In practice, the combination of uncertainty and variability can exhibit non-Gaussian features; however, the data-driven analysis from our prior work, \cite{dvorkin_2016f, lubin_2016}, suggests that the Gaussian assumption is sufficiently accurate for CCOPF computations over a large transmission grids. One can also improve the accuracy of the Gaussian representation by introducing uncertainty sets on its parameters, see \cite{lubin_2016}, or by using an out-of-sample analysis of \cite{bienstock_2014} to adjust parameters $\varepsilon_i^\uparrow$, $\varepsilon_i^\downarrow$,  $\varepsilon_{ij}^\uparrow$, and $\varepsilon_{ij}^\downarrow$ in a way that they would  match a given Gaussian distribution to a desired non-Gaussian distribution. The optimization Eq.~\eqref{CC-OPF}-\eqref{line_constr_hard2} assumes that  $p_l$ remains fixed, as in Eq.~\eqref{DC-OPF}-\eqref{line_constr_hard}, while  $p_g$ depends on  $\rho$. This dependency makes it possible for conventional generators to deviate from their original set points, $p_g^{(0)}$, to follow random wind power outputs according to a given control policy. This policy can be affine or nonlinear as discussed in Section \ref{subsec:affine} and  Section \ref{subsec:dead}, respectively. Eq.~\eqref{gen_constr_hard1} and \eqref{line_constr_hard2} are chance constrained equivalents of Eq.~\eqref{gen_constr_hard} and \eqref{line_constr_hard}, respectively. Parameters $\varepsilon_{i}^\uparrow$, $\varepsilon_{i}^\downarrow$, $\varepsilon_{ij}^\uparrow$, $\varepsilon_{ij}^\downarrow$ can be interpreted as  proxies for the fraction of time (probability) when a constraint violation can be tolerated, \cite{bienstock_2014}. Relative to \cite{bienstock_2014,lubin_2016}, the optimization Eq.~\eqref{CC-OPF}--\eqref{line_constr_hard2} contains a number of simplifications made for the sake of clarity. First, it assumes compulsory participation of conventional generators in the frequency control provision. Second, it ignores parameter uncertainty on the probability distribution characterizing wind power generation and correlation between distribution parameters at different nodes.

\subsection{Affine Frequency Control}
\label{subsec:affine}

The deterministic OPF and CCOPF formulations in Section~\ref{subsec:CCOPF} are stated for a quasi-stationary (balanced) system state, characterized by $\sum_{i \in \cal{V}} p_i^{(0)}=0$. Random fluctuations of wind power outputs, given by $\sum_{i \in \cal{V}_{\omega}} \rho_i \neq 0$, forces the system imbalance, i.e.  $\sum_{i \in \cal{V}} p_i^{(0)} \neq 0$. In other words, since this paper focuses on the optimization over time scales of minutes and longer,  sub-minute transient processes (evolution from pre-disturbed state to post-disturbed state) are ignored. This simplification makes it possible to treat the frequency and power unbalances as proportional quantities, so that if there is no frequency control, the system equilibrates within a few seconds at
\begin{equation}
\omega=\frac{\sum_{i\in{\cal V}_w} \rho_i}{\sum_{k\in{\cal V}_g}\gamma_k},
\label{omega-u}
\end{equation}
where $\omega$ is a deviation from the nominal frequency and $\gamma_k$ is the (natural) damping coefficient of the generator $k$.

When the affine primary frequency control is activated, the conventional generators respond as
\begin{equation}
\forall i\in{\cal V}_g:\quad p_i^{(0)}\to p_i=p_i^{(0)}-\alpha_i^{(1)} \omega_i.
\label{primary}
\end{equation}
This response of the conventional generators force the system to equilibrate (within the same transient time scales ranging from a few seconds to tens of seconds) at
\begin{equation}
\omega^{(1)}=\frac{\sum_{i\in{\cal V}_w} \rho_i}{\sum_{i\in{\cal V}_g}(\alpha_i^{(1)}+\gamma_i)},
\label{omega-1}
\end{equation}
where $\alpha^{(1)}_{i}$ is the primary droop coefficient (participation factor) of the conventional generator $i$. Following the primary frequency control, the secondary frequency control, also called Automatic Generation Control (AGC) \cite{jaleeli_1992}, would be deployed within a few minutes, thus resulting in an additional affine correction to  the power output of the conventional generators. Then, Eq.~\eqref{primary} is modified to account for the secondary frequency control deployment  as follows
\begin{equation}
\forall i\in{\cal V}_g:\quad p_i^{(0)}\to p_i=p_i^{(0)}-\alpha_i^{(1)}\omega^{(1)}-\alpha_i^{(2)} \omega^{(2)},
\label{secondary}
\end{equation}
where $\alpha^{(2)}_{i}$ is the secondary droop coefficient (participation factor) of the conventional generator $i$. Eq.~\eqref{secondary} assumes that the additional correction is distributed among conventional generators according to $\alpha_i^{(2)}$ and ensures that the system is globally balanced after both the primary and secondary frequency responses are fully deployed. Note that Eq.~\eqref{secondary} ignores the inter-area correction component to simplify notations. Combining  Eq.~\eqref{omega-1} and Eq.~\eqref{secondary} results in
\begin{equation}
\omega^{(2)}\sum\limits_{i\in{\cal V}_g} \alpha_i^{(2)}=\sum\limits_{i\in{\cal V}_w}\rho_i -\omega^{(1)}\sum\limits_{i\in{\cal V}_g}\alpha_i^{(1)}=\omega^{(1)}\sum\limits_{i\in{\cal V}_g}\gamma_i^{(1)},
\label{omega-2}
\end{equation}
It is noteworthy to note that the  secondary control should be considered as a centrally-controlled addition to the locally-managed primary control, see \cite{jaleeli_1992}.

Note that $\omega^{(1)}$ and $\omega^{(2)}$  can each be expressed in terms of $\sum_{i\in{\cal V}_w} \rho_i$  using Eq.~\eqref{omega-2}. These expressions can then be combined with Eq.~\eqref{secondary} to summarize Eq.~\eqref{flow_node}--\eqref{flow_line} 
\begin{equation}
\forall i\in{\cal V}: \sum\limits_{j\sim i} \beta_{ij}(\theta_i-\theta_j)=\left\{\begin{array}{cc}
p_i^{(0)}-\tilde{\alpha}_i \sum_{k\in{\cal V}_w}\rho_k, & i\in{\cal V}_g\\
\rho_i, & i\in{\cal V}_w\\ p_i,& i\in{\cal V}_l\end{array}\right.
\label{PF_fluct}
\end{equation}
where $\tilde{\alpha}_i$ stands for the renormalized droop coefficient of the conventional generators computed according to
\begin{equation}
\forall i\in{\cal V}_g:\quad \tilde{\alpha}_i=\frac{\alpha_i^{(1)}+\alpha_i^{(2)}\frac{\sum\limits_{k\in {\cal G}_g} \gamma_k^{(1)}}{\sum\limits_{l\in{\cal V}_g} \alpha_l^{(2)}}}{\sum\limits_{m\in{\cal V}_g}(\alpha_m^{(1)}+\gamma_m^{(1)})}.\label{tilde_alpha}
\end{equation}
The renormalized droop coefficients are subject to the following integrality constraint:
\begin{equation}
\sum_{i\in{\cal V}_g}\tilde{\alpha}_i=1.
\end{equation}

\subsection{Dead Zone in Chance-Constrained Primary Control}
\label{subsec:dead}

Eq.~\eqref{secondary}--\eqref{PF_fluct} assume that the primary control reacts to a frequency deviation of any size. Even though this assumption is applicable for some systems (e.g.  microgrids), it does not necessarily hold for large systems, where the primary frequency control is routinely kept untarnished during normal operations and is used only for quick and relatively rare response to large disturbances, e.g.  contingencies. In such systems, $\alpha_i^{(1)}$  depends on the size of the frequency deviation, $\omega^{(1)}$, and can be formalized as:
\begin{eqnarray}
\forall i\in {\cal V}_g:\quad \alpha_i^{(1)}\rightarrow \alpha_i^{(1)} \left\{\begin{array}{cc} 0,& |\omega^{(1)}|\leq \bar{\Omega} \\ 1,& \mbox{otherwise}\end{array}\right.,
\label{dead-zone}
\end{eqnarray}
where parameter $\overline{\Omega}$ is a frequency threshold (dead zone) for the primary frequency response. 
The value of this threshold can be manually chosen as it suits the needs of a particular system, \cite{dvorkin_2016}. Note that the dead zone makes the frequency control policy given by Eq.~\eqref{dead-zone} nonlinear; hence, one generally anticipates that ignoring the nonlinearity and using the affine policy would cause some inaccuracy.

To account for the dead zone in Eq.~\eqref{dead-zone}, we suggest the following modification of  \eqref{gen_constr_hard1}--\eqref{gen_constr_hard2}:

\begin{eqnarray}
&\forall i\in{\cal V}_g: &\!\mathbb{E}_{\rho}\left[\theta(\underline{p}_i-p_i)\theta(\bar{\Omega}-|\omega^{(1)}|)\right]
\leq \varepsilon_i^{(\downarrow,-)},\label{CC_gen_2a}\\
& &\!\mathbb{E}_{\rho}\left[\theta(\underline{p}_i-p_i)\theta(|\omega^{(1)}|-\bar{\Omega})\right]
\leq \varepsilon_i^{(\downarrow,+)},\label{CC_gen_2b}\\
& &\!\mathbb{E}_{\rho}\left[\theta(p_i-\overline{p}_i)\theta(\bar{\Omega}-|\omega^{(1)}|)\right]
\leq \varepsilon_i^{(\uparrow,-)},\label{CC_gen_2c}\\
& &\!\mathbb{E}_{\rho}\left[\theta(p_i-\overline{p}_i)\theta(|\omega^{(1)}|-\bar{\Omega})\right]
\leq \varepsilon_i^{(\uparrow,+)},\label{CC_gen_2d}
\end{eqnarray}
where $\theta(x)$ is the unit step function (also known as the Heaviside step function) such that $\theta(x)=1$, if $x>0$, and $\theta(x)=0$ otherwise. Eqs.~\eqref{CC_gen_2a}--\eqref{CC_gen_2d} incorporate both the nonlinear primary frequency response and the wind power generation statistics described by $\rho$, as given in Eq.~\eqref{wind_mean}-\eqref{wind_cov}. Note that the parameters $\varepsilon_i^{(\downarrow,-)}$, $\varepsilon_i^{(\downarrow,+)}$, $\varepsilon_i^{(\uparrow,-)}$,$\varepsilon_i^{(\uparrow,+)}$ are defined similarly to parameters  $\varepsilon_{i}^\uparrow$ and  $\varepsilon_{i}^\downarrow$ in Eq.~\eqref{gen_constr_hard1}-\eqref{gen_constr_hard2}.

Section \ref{sec:solution_approach} shows how Eq.~\eqref{CC_gen_2a}--\eqref{CC_gen_2b} can be represented in a computationally tractable form as one-dimensional integrals (with erf-functions in the integrands). The same transformation  is applicable for Eq.~\eqref{CC_gen_2c}--\eqref{CC_gen_2d}.

\subsection{Weighted CCOPF}
\label{subsec:weighted}

Following the method of \cite{roald_2015} we introduce weighted chance constraints imposed at the conventional generators:
\begin{eqnarray}
&\forall i\in{\cal V}_g: & \!\!\mathbb{E}_{\rho}\!\!\left[ \exp\!\left(\!-\frac{p_i}{\underline{p}_i}\!\right) \theta\left(\bar{\Omega}-|\omega^{(1)}|\right)\right]\!\leq\! \varepsilon_i^{(\downarrow,-)},
\label{WCC_gen_2a}\\
&& \!\!\mathbb{E}_{\rho} \!\!\left[ \exp\!\left(\!-\frac{p_i}{\underline{p}_i}\!\right) \theta\left(|\omega^{(1)}|-\bar{\Omega}\right)\right]\!\leq\! \varepsilon_i^{(\downarrow,+)},
\label{WCC_gen_2b}\\
&& \!\!\mathbb{E}_{\rho} \!\!\left[ \exp\!\left(\!\frac{p_i}{\overline{p}_i}\!\right) \theta\left(\bar{\Omega}-|\omega^{(1)}|\right)\right]\!\leq\! \varepsilon_i^{(\uparrow,-)},
\label{WCC_gen_2c}\\
&& \!\!\mathbb{E}_{\rho} \!\! \left[ \exp\!\left(\!\frac{p_i}{\overline{p}_i}\!\right) \theta\left(|\omega^{(1)}|-\bar{\Omega}\right)\right]\!\leq\! \varepsilon_i^{(\uparrow,+)}.
\label{WCC_gen_2d}
\end{eqnarray}
The weighted chance constraints are advantageous for the following three reasons. First, this form allows the freedom in differentiating effects of small and large violations \cite{roald_2015}. Second, the resulting constraints are convex regardless of the input statistics (Gaussian or not) \cite{roald_2015}. Third,  it offers a computational advantage as the expectations on the left-hand side of Eqs.~\eqref{WCC_gen_2a}-\eqref{WCC_gen_2d} are stated explicitly in terms of the (well tabulated) erf-functions.

Then our CCOPF-PFR formulation with the weighted chance constraints is as follows:
\begin{eqnarray}
\hspace{-5mm}&&\min\limits_{p_g^{(0)}} \sum_{i\in{\cal V}_g} \mathbb{E}_{\rho}\left[C_i(p_i)\right]
\label{CCOPF-PFR}\\
\hspace{-5mm}&  \forall i\in{\cal V}_g: &\mbox{Eqs.~\eqref{WCC_gen_2a}--\eqref{WCC_gen_2d}}  \label{CCOPF-PFR2} \\
\hspace{-5mm} & \forall \{i,j\}\in{\cal E}: &  \hspace{-3mm} \mathbb{E}_{\rho} \hspace{-1mm}\left[ \exp\left(\!-\frac{\phi_{ij}}{\underline{\phi}_{ij}}\right) \theta\left(\bar{\Omega}\!-\!|\omega^{(1)}|\!\right)\right] \hspace{-2mm} \leq \!\varepsilon_{ij}^{(\downarrow,-)}, \label{WCC_line_2a}\\
\hspace{-3mm}&& \hspace{-3mm} \mathbb{E}_{\rho}\hspace{-1mm}\left[ \exp\left(\!-\frac{\phi_{ij}}{\underline{\phi}_{ij}}\right) \theta\left(|\omega^{(1)}|\!-\!\bar{\Omega}\!\right)\right]\hspace{-2mm}\leq \!\varepsilon_{ij}^{(\downarrow,+)},
\label{WCC_line_2b}
\end{eqnarray}
\begin{eqnarray}
&& \mathbb{E}_{\rho}\left[ \exp\left(\frac{\phi_{ij}}{\overline{\phi}_{ij}}\right) \theta\left(\bar{\Omega}-|\omega^{(1)}|\right)\right]\leq \varepsilon_{ij}^{(\uparrow,-)},
 \label{WCC_line_2c}\\
&& \mathbb{E}_{\rho}\left[ \exp\left(\frac{\phi_{ij}}{\overline{\phi}_{ij}}\right) \theta\left(|\omega^{(1)}|-\bar{\Omega}\right)\right]\leq \varepsilon_{ij}^{(\uparrow,+)}.
\label{CCOPF-PFR_end} 
\end{eqnarray}

Note that parameters $\varepsilon_{ij}^{(\downarrow,-)}$, $\varepsilon_{ij}^{(\downarrow,+)}$, $\varepsilon_{ij}^{(\uparrow,-)}$,$\varepsilon_{ij}^{(\uparrow,+)}$ are defined similarly to parameters  $\varepsilon_{i}^\uparrow$ and  $\varepsilon_{i}^\downarrow$ in Eqs.~\eqref{gen_constr_hard1}-\eqref{gen_constr_hard2}.

\section{Solution Approach} \label{sec:solution_approach}

This Section describes how the weighted chance constraints \eqref{CCOPF-PFR2}-\eqref{CCOPF-PFR_end} can be computed efficiently. The process is illustrated on Eqs.~\eqref{WCC_gen_2a}-\eqref{WCC_gen_2b} and can be extended to other chance constraints. The left-hand sides of Eq.~\eqref{WCC_gen_2a}--\eqref{WCC_gen_2b} depend on $\rho_i$ and $\overline{\Omega}$ attaining non-zero values if  $\bar{\Omega}>|\omega^{(1)}|$ and $\bar{\Omega}<|\omega^{(1)}|$, respectively. Thus, the left-hand sides can be restated as expectations over two distinct Gaussian distributions defined by the following means and covariances:
\begin{eqnarray}
&& \mathbb{E}\left[\omega^{(1)}\right]_\sigma\doteq \Omega_\sigma=\frac{\sum_{i\in{\cal V}_w}\bar{\rho}_i}{\sum_{k\in{\cal V}_g}(\sigma\alpha^{(1)}_k+\gamma_k)}, \label{omega_mean}\\
&& \mathbb{E}\left[\left(\omega^{(1)}-\Omega_\sigma\right)^2\right]_\sigma\doteq \Theta^{(\omega,\omega)}_\sigma\nonumber\\ && =\frac{\sum_{i\in{\cal V}_w}R_i}{(\sum_{k\in{\cal V}_g}(\sigma\alpha^{(1)}_k+\gamma_k))^2}, \label{omega_cov}\\
&&\forall i\in{\cal V}_g: \nonumber\\
&& \mathbb{E}\left[p_i\right]_\sigma\doteq P_{\sigma;i}=p^{(0)}_i-\tilde{\alpha}_{i;\sigma} \sum_{j\in{\cal V}_w}\bar{\rho}_j, \label{p_mean_sigma}\\
&& \mathbb{E}\left[\left(p_i-P_{\sigma;i}\right)\left(\omega^{(1)}-\Omega_\sigma\right)\right]_\sigma\doteq \Theta^{(p_i,\omega)}_\sigma\nonumber\\ && =-\frac{\tilde{\alpha}_{i;\sigma}\sum_{j\in{\cal V}_w} R_j}{\sum_{k\in{\cal V}_g}(\sigma \alpha_k^{(1)}+\gamma_k)}
\label{omega-p_cov}\\
&& \mathbb{E}\left[\left(p_i-P_{\sigma;i}\right)^2\right]_\sigma\doteq \Theta^{(p_i,p_i)}_\sigma= \left(\tilde{\alpha}_{i;\sigma}\right)^2\sum_{j\in{\cal V}_w} R_j,
\label{p-p_cov}
\end{eqnarray}
where  $\sigma=\{0,1\}$ distinguish the case ``$\omega^{(1)}$-in-range" and the case ``$\omega^{(1)}$-off-range". If $\sigma=0$,  $\tilde{\alpha}_{i;\sigma}=\tilde{\alpha}_{i}$. If $\sigma=1$, $\tilde{\alpha}_{i}$ is derived from \eqref{tilde_alpha} using the replacement $\alpha^{(1)}=0$. It is also useful to introduce the so-called precision matrices defined as the inverse $(2\times 2)$ matrices of the covariance matrices
\begin{eqnarray}
&\forall i\in{\cal V}_g:& \Phi_{\sigma;i}= \left(\begin{array}{cc} \Phi^{(p_i,p_i)}_\sigma & \Phi^{(p_i,\omega)}_\sigma\\
\Phi^{(p_i,\omega)}_\sigma & \Phi^{(\omega,\omega)}_\sigma\end{array}\right)\doteq (\Theta_{\sigma;i})^{-1}\nonumber\\
&&=\frac{\left(\begin{array}{cc}  \Theta^{(\omega,\omega)}_\sigma & -\Theta^{(p_i,\omega)}_\sigma\\
-\Theta^{(p_i,\omega)}_\sigma & \Theta^{(p_i,p_i)}_\sigma \end{array}\right)}{\Theta^{(p_i,p_i)}_\sigma \Theta^{(\omega,\omega)}_\sigma - (\Theta^{(p_i,\omega)}_\sigma)^2}.
\label{prec_matrix}
\end{eqnarray}
Next, Eqs.~\eqref{WCC_gen_2a} can be simplified as:

\begin{eqnarray}
&& \int\limits_{-\infty}^{+\infty}dx \left(\exp\left(\frac{x}{\underline{p}_i-P_{0;i}}\right)-1\right) \times \label{WCC_gen_2a_bis}\\ \nonumber
\end{eqnarray}
\begin{eqnarray}
&& \times \int\limits_{-\bar{\Omega}-\Omega_0}^{\bar{\Omega}-\Omega_0} dy \frac{\sqrt{\mbox{det}(\Phi_{0;i})}\exp\left(-\frac{1}{2}(x,y)\Phi_{0;i}\left(\begin{array}{c} x\\ y\end{array}\right)\right)}{2\pi}.\nonumber
\end{eqnarray}
Eq.~\eqref{WCC_gen_2a_bis} is convex with respect to the set points of the conventional generators. Note that chance constraints \eqref{WCC_gen_2b}--\eqref{WCC_gen_2d}, \eqref{WCC_line_2a}--\eqref{CCOPF-PFR_end} can be converted to similar convex expressions.

\section{Case Study} \label{sec:case_study}

The proposed CCOPF-PFR is compared to the standard CCOPF over a modification of the 118-bus Reliability Test System \cite{ieee118}. The test system includes 54 conventional generation resources and 186 transmission lines. Additionally, this system includes 9 wind farms with the forecasted total power output of 1,053 MW as itemized in Table~\ref{tab:wind}.   The mean and standard deviation of the wind power  outputs are set to 0\% and 10\% of the power forecast at each wind farm. The power flow limit of each transmission line is reduced by 25\% of its rated value and the active power demand is increased by 10\% at each bus. The droop coefficients (participation factors) of each conventional generator for each regulation interval are set to $1/N_g$, where $N_g$ is the number of conventional generators, i.e. $\alpha_i^{(1)}=\alpha_i^{(2)} = 1/N_g, \forall i \in \mathcal{V}_g$. The value of the dead zone for the primary frequency response is 100 MW and the likelihood of the constraint violations is set to  $\varepsilon_i^\downarrow=\varepsilon_i^\uparrow=\varepsilon, \forall i \in {\cal V}_g, $ and $\varepsilon_{ij}^\uparrow=\varepsilon_{ij}^\downarrow=\varepsilon, \forall \left\{i,j\right\} \in {\cal E} $. Both the CCOPF and CCOPF-PFR formulations are implemented in Julia \cite{miles_julia2015} using the JumpChance package resolved using a 1.6 Ghz Intel Core i5 processor with 8GB of RAM.

\begin{center}
\begin{table}[h]
 \captionsetup{justification=centering, labelsep=period, font=footnotesize, textfont=sc}
\caption{Wind Power Forecast at  Wind Farms  ($\overline{\omega}$, MW)}
\begin{tabular}{ c| c c c c c c c c c c c }
\hline
\hline
bus \#  &3 &8 &11 & 20& 24 &38& 43 & 49 & 50\\
 \hline
$\overline{\omega}$  & 70& 147 &102 & 105&113 & 250&118 & 76 & 72\\
\hline
\hline
\end{tabular}
\label{tab:wind}
\end{table}
\vspace{-10mm}
\end{center}

\subsection{Comparison of the CCOPF and CCOPF-PFR solution}
\vspace{-7mm}
\begin{center}
\begin{table}[b]
 \captionsetup{justification=centering, labelsep=period, font=footnotesize, textfont=sc}
\caption{Objective Function Values and CPU Times for the CCOPF and CCOPF-PFR Formulations}
\begin{tabular}{ c| c c| c c }
\hline
\hline
\multirow{2}{*}{$\varepsilon$} & \multicolumn{2}{c|}{Objective function, \$} &
\multicolumn{2}{c}{CPU time, s}\\
\cline{2-5}
 &CCOPF & CCOPF-PFR &CCOPF  & CCOPF-PFR\\
 \hline
$10^{-1}$ &91,504.8 & 94,757.3 (+3.55\%)$^*$ &10.4  &17.3\\
$10^{-2}$  &92,984.6 & 94,861.5 (+2.02\%)$^*$ & 11.9 & 18.1\\
$10^{-3}$   &94658.7 & 96,045.3 (+1.46\%)$^*$ &12.1 & 19.4\\
$10^{-4}$       &97,615.1 & 98,101.4 (+0.49\%)$^*$  &21.7  & 18.9\\
\hline
\hline
\multicolumn{5}{l}{$^*$ -- the percentage values are relative to the CCOPF formulation.}
\end{tabular}
\label{tab:cost_time}
\end{table}
\end{center}
First, the CCOPF and CCOPF-PFR formulations are solved for different values of $\varepsilon$. Table~\ref{tab:cost_time} compares the objective functions of these formulations and their CPU times. In general, the CCOPF-PFR consistently yields a more expensive solution since it has a more constrained feasible region due to the active PFR constraints. As the value of parameter $\varepsilon$ reduces so does the relative difference between the objective function of the CCOPF and CCOPF-PFR formulations. This observation suggests that a nonlinear primary frequency response policy comes at a lower cost for risk-averse OPF solutions. At the same time, the computing times reported in Table~\ref{tab:cost_time}  suggest that a higher level of modeling accuracy in the  CCOPF-PFR formulation comes at a minimal increase in the computation time.

Second, the CCOPF and CCOPF-PFR solutions presented in Table~\ref{tab:cost_time} are tested over a set of 10,000 random realizations of wind power forecast errors, where each random realization is sampled using the multivariate Gaussian distribution. In each test, the decisions produced by the CCOPF and CCOPF-PFR  are fixed and the constraint violations are calculated.

\begin{figure}[t]
    \centering
    \includegraphics[width=\columnwidth]{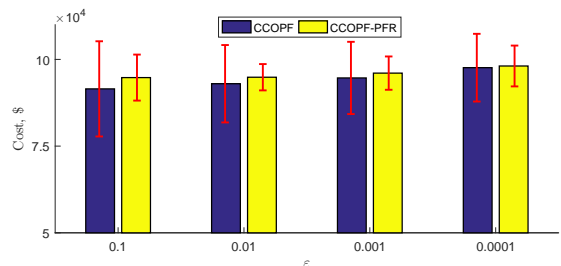}
    \caption{\small Comparison of the CCOPF and CCOPF-PFR formulations in terms of the expected costs (vertical bars) and standard deviations (error bars) for different values of parameter $\varepsilon$.  }
    \label{fig:cost}
\vspace{-5mm}
\end{figure}

\begin{figure}[t]
    \centering
    \includegraphics[width=\columnwidth]{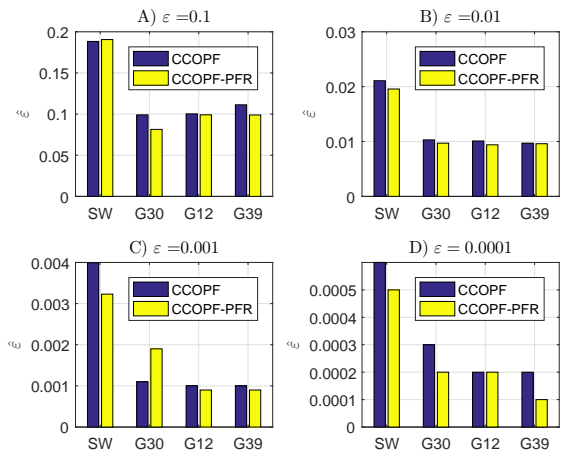}
    \caption{\small  Empirical violations of chance constraints on conventional generators, $\hat{\varepsilon}$, for different values of parameter $\varepsilon$. Label SW (system-wide) denotes the fraction of realizations that lead to a violation of at least one chance constraint in the entire transmission system. Labels G30, G12 and G38 denote the fraction of realizations that lead to violations of chance constraints for the maximum power output at  conventional generators  G30 ($\overline{p}_i$ = 805.2 MW), G12 ($\overline{p}_i$ = 413.9 MW) and G38 ($\overline{p}_i$ =104.0 MW)}
    \label{fig:viol}
    \vspace{-7mm}
\end{figure}

Fig.~\ref{fig:cost} displays the excepted operating cost and its standard deviation observed over 10,000 tests for the CCOPF and CCOPF-PFR formulations. For both formulations the expected operating cost and its standard deviation monotonically change with the value of parameter $\varepsilon$. Thus, the expected operating cost under both formulations gradually increases for higher values of parameter $\varepsilon$, while the standard deviation reduces. Notably, in all instances displayed in Fig.~\ref{fig:cost} the expected cost of the CCOPF-PFR formulation is greater than the expected cost of the CCOPF formulation. As in the cost results presented in Table~\ref{tab:cost_time}, the gap between the expected costs of both formulations reduces for higher values of parameter $\varepsilon$. On the other hand, the CCOPF-PFR formulation leads to a lower standard deviation in all instances, which suggests that the CCOPF-PFR formulation is more robust and cost-efficient for accommodating relatively large  wind power forecast errors. A more expensive and conservative CCOPF-PFR solution leads to less violations of chance constraints on conventional generators as shown in Fig.~\ref{fig:viol}. The number of violations reduces for the entire system and for individual generators of different sizes. Therefore, the CCOPF-PFR formulation is more effective in accommodating deviations from the forecasted values.  Furthermore, there is no noticeable difference in violations of the chance constraints over line flows between the CCOPF and CCOPF-PFR formulation. This observation suggests that the main effect of  CCOPF-PFR is in improving compliance on the supply side with the operating limits.

\section{Conclusion} \label{sec:conclusion}

In this paper, the CCOPF formulation from  \cite{bienstock_2014} has been enhanced to  explicitly model constraints associated with the primary frequency response based on a nonlinear speed-droop governor characteristic of conventional generators. The proposed CCOPF-PFR formulation has been compared to the original CCOPF formulation on a modification of the 118-bus IEEE Reliability Test System \cite{ieee118}. This comparison indicates that modeling a nonlinear speed-droop governor characteristic leads to only a rather modest increase of the expected operating cost, while improving adaptability of the dispatch solutions to relative large deviations from the forecast. The increased adaptability of the CCOPF-PFR formulation is observed in reduction of the chance constraints violations on the conventional generators and it is also seen in lower standard deviations of the operating cost.  We have also observed that the  proposed CCOPF-PFR and standard CCOPF formulations are comparable in terms of  required computational resources.

This work can be extended in several ways:
\begin{itemize}
    \item The PFR constraints assume that generators instantly react to a power imbalance, i.e. there is no time delay, which can be observed in practice \cite{jaleeli_1992}. Modeling this delay is a possible extension of the proposed work aimed at improving the accuracy. 

    \item PFR constraints can be generalized to explicitly account for instant power flow fluctuations over transmission lines adjusted to the generator. This can be without additional communication constraints by using local measurements.
        
    \item The proposed CCOPF-PFR model can be enhanced to include an endogenous contingency reserve assessment, e.g. a probabilistic security-constrained framework \cite{dvorkin_2017d}. The proposed primary frequency response constraints can be used to  accurately estimate the minimum response requirement and its allocation instead of using the deterministic heuristics \cite{dvorkin_2016}.
        
        \item The proposed CCOPF-PFR model relies on lossless DC power flows, which needs to and seemingly can be extended to account for power losses, reactive power  flows, and voltage fluctuations via linear or quadratic AC power flow approximations.
            
\end{itemize}

\bibliographystyle{IEEEtran}
\bibliography{ConferenceBib}

\end{document}